\documentclass{article}
\usepackage{amssymb}
\newtheorem{theorem}{Theorem}
\newtheorem{lemma}[theorem]{Lemma}

\newtheorem{definition}[theorem]{Definition}
\newtheorem{example}[theorem]{Example}
\newtheorem{remark}[theorem]{Remark}

\def\QED{\quad\blackslug\lower 8.5pt\null}

\newcommand{\crazy}[2]{\displaystyle{\mathop{#1}_{#2}}
\vphantom{\displaystyle{#1}}}

\begin{document}

\begin{center}
{\Large\bf GOURSAT'S {\mathversion{bold} $(n+1)$}-WEBS }

\vspace*{4mm}

{\large VLADISLAV V. GOLDBERG}
\end{center}

{\small{\bf Abstract}.
{\em  We consider the
 Goursat's $(n+1)$-webs of codimension one  of two kinds
 on  an $n$-dimensional manifold. They are
 characterized by the specific  closed form equations
 or by two special relations between components
 of the torsion tensor of the web. These relations
  allow us to establish a connection with solutions of
 two systems of nonlinear second-order PDEs
 investigated by Goursat in 1899.
  The integrability conditions
 of some distributions invariantly associated
 with both kinds of Goursat's $(n+1)$-webs are
 also investigated.}

 \vspace*{8mm}

\section{The principal equations of  $(n+1)$-webs}

\textbf{1.} Let $W(n+1,n,1)$ be an $(n+1)$-web defined on a
differentiable manifold $X^{n}$ of dimension $n$
by $n + 1$ foliations $\lambda_{\xi},\;\; \xi = 1,...,n+1,$
of codimension one. Each foliation $\lambda_{\xi}$
can be defined by the completely integrable system of Pfaffian
equations
\begin{equation}\label{eq:1}
 \crazy{\omega}{\xi}=0,\;\; \xi=1,\ldots,n+1.
\end{equation}
The 1-forms $\crazy{\omega}{\alpha},\;\;
\alpha = 1, \ldots, n$,
define a co-frame in the tangent bundle $T(X^{n})$ and satisfy
the following structure equations:
\begin{equation}\label{eq:2}
d\crazy{\omega}{\alpha} = \crazy{\omega}{\alpha} \wedge
\omega +  \sum_{\beta\neq\alpha} \;
\crazy{a}{\alpha\beta} \;
\crazy{\omega}{\alpha} \wedge \crazy{\omega}{\beta},
\end{equation}
\begin{equation}\label{eq:3}
d\omega = \sum_{\alpha,\beta = 1}^{n}
\crazy{b}{\alpha\beta}\;\crazy{\omega}{\alpha}
\wedge \crazy{\omega}{\beta},
\end{equation}
where the quantities $\crazy{a}{\alpha\beta}$ and
$\crazy{b}{\alpha\beta}$ are connected by
certain relations (see [G 73] or [G 74] or [G 88], Section \textbf{1.2}).
We indicate some of these relations:
\begin{equation}\label{eq:4}
\nabla\crazy{a}{\alpha\beta}= \sum_{\gamma = 1}^{n}
(\crazy{a}{\alpha\beta\gamma} +
\crazy{a}{\alpha\beta} \;
\crazy{a}{\gamma\alpha} +
\crazy{a}{\alpha\beta} \;
\crazy{a}{\beta\gamma}) \;\crazy{\omega}{\gamma} ,
\;\; \alpha \neq \beta,
\end{equation}
\begin{equation}\label{eq:5}
\crazy{a}{\alpha\beta} =
\crazy{a}{\beta\alpha},
\end{equation}
\begin{equation}\label{eq:6}
\sum_{\alpha,\beta=1}^{n} \crazy{a}{\alpha\beta} = 0,
\end{equation}
\begin{equation}\label{eq:7}
\crazy{b}{\alpha\beta} =
\frac{1}{2}(\crazy{a}{\gamma\alpha\beta} -
\crazy{a}{\beta\gamma\alpha}),\;\;
\gamma \neq \alpha,\beta,
\end{equation}
\begin{equation}\label{eq:8}
\crazy{b}{\alpha\beta} = \crazy{b}{\alpha\gamma} + \crazy{b}{\gamma\beta},
\end{equation}
where $\nabla\crazy{a}{\alpha\beta} =
d\crazy{a}{\alpha\beta}
- \crazy{a}{\alpha\beta} \; \omega.$ By (7), we have
\begin{equation}\label{eq:9}
\crazy{b}{\alpha\beta} = -\crazy{b}{\beta\alpha}.
\end{equation}

The quantities $\crazy{a}{\alpha\beta}$ and
$\crazy{b}{\alpha\beta}$ form tensor fields in
the tangent bundle
 $T(X^n)$ which are called respectively the
{\it torsion} and {\it curvature tensors} of the
web $W(n+1,n,1)$.

\section{Goursat's $(n+1)$-webs of first kind}

\textbf{2.} Suppose that $n \geq 4$.
 Consider two pairs of three-dimensional distributions defined
by the following Pfaffian equations:
\begin{equation}\label{eq:10}
  \crazy{a}{13}\; \crazy{\omega}{3} + \crazy{a}{14} \;\crazy{\omega}{4}
  =0, \; \crazy{\omega}{\sigma} = 0,
\end{equation}
\begin{equation}\label{eq:11}
  \crazy{a}{23}\; \crazy{\omega}{3} + \crazy{a}{24}\; \crazy{\omega}{4}
  =0, \; \crazy{\omega}{\sigma} = 0,
\end{equation}
and
\begin{equation}\label{eq:12}
  \crazy{a}{13} \;\crazy{\omega}{1} + \crazy{a}{23}\; \crazy{\omega}{2}
  =0, \; \crazy{\omega}{\sigma} = 0,
\end{equation}
\begin{equation}\label{eq:13}
  \crazy{a}{14}\; \crazy{\omega}{1} + \crazy{a}{24}\; \crazy{\omega}{2}
  =0, \; \crazy{\omega}{\sigma} = 0,
\end{equation}
where $\sigma = 5, \ldots , n.$

In general, the distributions defined by each of equations (10)--(13)
are not integrable.

It is easy to see that
{\em the distribution $(10)$ and $(11)$
$($or $(12)$ and $(13))$ coincide
if and only if the torsion tensor
$\crazy{a}{\alpha\beta}$ of the web $W (n + 1, n, r)$
satisfies the following condition}:
\begin{equation}\label{eq:14}
\crazy{a}{13} \; \crazy{a}{24} - \crazy{a}{14} \; \crazy{a}{23} = 0.
\end{equation}

One can also see that {\em condition $(14)$ is necessary
and sufficient for the distribution defined
by the system of equations $(10)$ and  $(11)$ $($or $(12)$ and $(13))$
to be three-dimensional} (for a general web $W$ these
distributions are two-dimensional).

We will call the webs $W(n+1,n,1)$ satisfying
condition (14) the {\em Goursat webs of the first kind}.
The reason for this name will be clear later.

\begin{lemma} The Pfaffian derivatives of the torsion tensor
$\crazy{a}{\alpha\beta}$ of the Goursat web $W (n + 1, n, r)$
of the first kind and the components of the torsion tensor
$\crazy{a}{\alpha\beta}$ itself satisfy the following relations:
\begin{equation}\label{eq:15}
\crazy{a}{24} \;\crazy{a}{13c} + \crazy{a}{13}\; \crazy{a}{24c}
 - \crazy{a}{14} \;\crazy{a}{23c} - \crazy{a}{23} \;\crazy{a}{14c} = 0,
 \;\; c = 1, 2, 3, 4.
\end{equation}
\end{lemma}

{\sf Proof.} In fact, differentiating condition (14), we find that
\begin{equation}\label{eq:16}
\crazy{a}{24}\; \nabla \crazy{a}{13} +  \crazy{a}{13}
\; \nabla \crazy{a}{24} - \crazy{a}{14} \; \nabla \crazy{a}{23}
- \crazy{a}{23} \; \nabla \crazy{a}{14} = 0.
\end{equation}
By (4), (10), and (14) and  by linear independence
of the forms $\crazy{\omega}{c}, \, c = 1, 2, 3, 4$,
equation (16) implies  conditions (15). \rule{3mm}{3mm}

\begin{theorem}
For the  Goursat webs $W (n + 1, n, 1)$
of the first kind
  the two-dimensional distribution defined
by equations $(10)$ and $(12)$  $($or $(11)$ and $(13))$ is integrable.
\end{theorem}

{\sf Proof.} In fact, denote by $\theta$ and $\rho$ the left-hand
sides of equations (10) and (12):
$$
\theta =  \crazy{a}{13} \; \crazy{\omega}{3} + \crazy{a}{14}
\;\crazy{\omega}{4}, \;\; \rho =
 \crazy{a}{13} \;\crazy{\omega}{1} + \crazy{a}{23}
\; \crazy{\omega}{2}.
$$

Since  $W (n + 1, n, 1)$ is a Goursat web of the first kind,
the components of its torsion tensor and their Pfaffian
derivatives satisfy conditions (14) and (15). Applying these
 conditions, one can easily prove that
 $$
 \left\{
 \begin{array}{ll}
 d \theta \equiv 0 \pmod{\theta, \rho}, \\
  \\
 d \rho \equiv 0 \pmod{\theta, \rho}.
 \end{array}
 \right.
 $$
Thus, the two-dimensional distribution defined
by the equations $\theta = 0, \; \rho = 0$ is integrable.
\rule{3mm}{3mm}

\textbf{3.} Suppose that in some domain $D \subset  X^{n}$,
a web $W (n + 1, n, 1)$ is defined by
the closed form  equations
\begin{equation}\label{eq:17}
x_{n+1} = F (x_1,\ldots, x_n), \;\; \mbox{det}\left(
\frac{\partial F}{\partial x_{\alpha}} \right) \neq 0.
\end{equation}
It is proved in [G 76] (see also [G 88], Section {\bf 4.1}) that
\begin{equation}\label{eq:18}
\crazy{a}{\alpha\beta} = \frac{F_{\alpha\beta}}{F_\alpha F_\beta},
\end{equation}
and that by (18),  conditions (14) are equivalent to the following
second-order nonlinear partial
differential equation:
\begin{equation}\label{eq:19}
\frac{F_{13}}{F_{14}} = \frac{F_{14}}{F_{24}}.
\end{equation}

Goursat [Go 99] considered such an equation for $n = 4, 5$.
This is the reason we named  webs satisfying condition
(14)  Goursat webs.

The following theorem follows from [Go 99].

\begin{theorem}
 For the  Goursat web $W (n + 1, n, 1)$ of the first kind,
 closed form equation $(17)$ takes the form
\begin{equation}\label{eq:20}
\left\{
\begin{array}{ll}
x_{n+1} = \phi (x_1, x_2, x_5, \ldots, x_n, a)
      +  \psi (x_3, x_4, x_5, \ldots, x_n, a),\\
      \\
    \displaystyle  \frac{\partial \phi}{\partial a} +
\frac{\partial \psi}{\partial a} = 0,
\end{array}
\right.
\end{equation}
where $\phi$ and $\psi$ are arbitrary functions of $n-1$ variables each
satisfying the second equation.
\end{theorem}

{\sf Proof.} In fact, it is proved in [Go 99] that the general
solution of equation (19) has the form (20). The only difference
between [Go 99] and our considerations is that in [Go 99]
$n = 4$ and $n = 5$ while we consider the general case
$n \geq 4$. \rule{3mm}{3mm}

\section{Goursat's $(n+1)$-webs of second kind}

\textbf{4.} Suppose that $n \geq 5$. Consider the distribution
$\Delta_2$ defined by the following equations:
\begin{equation}\label{eq:21}
\left\{
\begin{array}{ll}
 \crazy{\omega}{3} + \crazy{\omega}{4} +  \crazy{\omega}{5}= 0,
  \\
  \crazy{a}{13}\; \crazy{\omega}{3} + \crazy{a}{14} \;\crazy{\omega}{4}
   + \crazy{a}{15}\; \crazy{\omega}{5}= 0,
  \\
  \crazy{a}{23}\; \crazy{\omega}{3} + \crazy{a}{24}\; \crazy{\omega}{4}
+ \crazy{a}{25}\; \crazy{\omega}{5} =0, \\
  \crazy{\omega}{\sigma} = 0, \;\; \sigma = 6, \ldots, n.
\end{array}
\right.
\end{equation}
Note that equation (1) and the first equation of (21) implies that
\begin{equation}\label{eq:22}
 \crazy{\omega}{1} +  \crazy{\omega}{2} = 0.
\end{equation}
Equations (21) and (22) imply that,  in general, the
distribution $\Delta_2$ is one-dimensional.

Consider also the distribution $\Delta_3$
defined by the equations
\begin{equation}\label{eq:23}
\left\{
\begin{array}{ll}
 \crazy{\omega}{3} + \crazy{a}{31}\;\crazy{\omega}{1} + \crazy{a}{32}\; \crazy{\omega}{2}= 0,
  \\
  \crazy{\omega}{3} + \crazy{a}{41} \;\crazy{\omega}{1}
   + \crazy{a}{42}\; \crazy{\omega}{2}= 0,
  \\
   \crazy{\omega}{3} + \crazy{a}{41}\; \crazy{\omega}{1}
+ \crazy{a}{42}\; \crazy{\omega}{2} =0, \\
  \crazy{\omega}{\sigma} = 0, \;\; \sigma = 6, \ldots, n.
\end{array}
\right.
\end{equation}
In general, the distribution $\Delta_3$ is
two-dimensional, and
the distributions defined by equations (21) and (23)
are not integrable.

Note that the $3 \times 3$ matrices of coefficients of the first
three equations of $(21)$ and $(23)$ are transposes of each other.
It is easy to see that {\em the distributions $\Delta_2$ and $\Delta_3$
defined by equations $(21)$ and $(23)$
are two- and three-dimensional, respectively, if and only if}
\begin{equation}\label{eq:24}
\left|
\begin{array}{lll}
1 & 1 & 1\\
\crazy{a}{13} & \crazy{a}{14} & \crazy{a}{15} \\
\crazy{a}{23} & \crazy{a}{24} & \crazy{a}{25}
\end{array}
\right| = 0.
\end{equation}
This explains why we used the notation $\Delta_2$ and $\Delta_3$
for these two distributions: the lower index indicates the
distribution dimension.

We will call the webs $W(n+1,n,1)$ satisfying
 condition (24)
 the {\em Goursat webs of the second kind}.
The reason for this name will be clear later.

It is easy to see that condition (24) is equivalent to
the condition
\begin{equation}\label{eq:25}
\left|
\begin{array}{ll}
\crazy{a}{13} & \crazy{a}{14}  \\
\crazy{a}{23} & \crazy{a}{24}
\end{array}
\right| +
\left|
\begin{array}{ll}
\crazy{a}{14} & \crazy{a}{15}  \\
\crazy{a}{24} & \crazy{a}{25}
\end{array}
\right| +
\left|
\begin{array}{ll}
\crazy{a}{15} & \crazy{a}{13}  \\
\crazy{a}{25} & \crazy{a}{23}
\end{array}
\right|
= 0,
\end{equation}
or to the condition
\begin{equation}\label{eq:26}
\crazy{a}{13} (\crazy{a}{24} - \crazy{a}{25})
+ \crazy{a}{14} (\crazy{a}{25} - \crazy{a}{23})
+ \crazy{a}{15} (\crazy{a}{23} - \crazy{a}{24})
+ \crazy{a}{23} (\crazy{a}{15} - \crazy{a}{14})
+ \crazy{a}{24} (\crazy{a}{13} - \crazy{a}{15})
+ \crazy{a}{25} (\crazy{a}{14} - \crazy{a}{13}) = 0,
\end{equation}
or to the condition
\begin{equation}\label{eq:27}
\displaystyle \frac{\crazy{a}{pa} - \crazy{a}{pb}}{\crazy{a}{pa} -
\crazy{a}{pc}} =
\frac{\crazy{a}{qa} - \crazy{a}{qb}}{\crazy{a}{qa} -
\crazy{a}{qc}}, \;\; p, q = 1, 2;\; a, b, c = 3, 4, 5; \; p \neq q;
\;a \neq b, c; \; b \neq c.
\end{equation}

\begin{theorem}
 For the  Goursat web $W (n + 1, n, 1)$ of the second kind,
 closed form equation $(17)$ take the form
\begin{equation}\label{eq:28}
\left\{
\begin{array}{ll}
x_{n+1} = \phi (x_1, x_2, x_6, \ldots, x_n, a,
        \psi (x_3, x_4, x_5, x_6, \ldots, x_n, a)),\\
      \\
    \displaystyle  \frac{\partial \phi}{\partial a} +
\frac{\partial \phi}{\partial \psi} \,\frac{\partial \psi}{\partial a} = 0,
\end{array}
\right.
\end{equation}
where $\phi$ and $\psi$ are arbitrary functions of
$n-1$ variables each, and $a$ is an arbitrary parameter.
\end{theorem}

{\sf Proof.} In fact, by  (18),
conditions (24) are equivalent to the following
second-order nonlinear partial
differential equation:
\begin{equation}\label{eq:29}
\left|
\begin{array}{lll}
F_3 & F_4 & F_5\\
F_{13} & F_{14} & F_{15} \\
F_{23} & F_{24} & F_{25}
\end{array}
\right| = 0.
\end{equation}
It is proved in [Go 99] that the general
solution of equation (29) has the form (28). The only difference
between [Go 99] and our considerations is that in [Go 99]
$n = 5$ while we consider the general case
$n \geq 5$. \rule{3mm}{3mm}

\textbf{5.} Let us find the differential consequences
of conditions (27). First, we will prove the following lemma.

\begin{lemma}
For the  Goursat web $W (n + 1, n, 1)$ of the second kind,
the following identities hold:
\begin{equation}\label{eq:30}
\crazy{a}{pa} (\crazy{a}{qc}  - \crazy{a}{qb}) +
\crazy{a}{pb} (\crazy{a}{qa}  - \crazy{a}{qc}) +
\crazy{a}{pc} (\crazy{a}{qb}  - \crazy{a}{qa}) = 0,
\end{equation}
\begin{equation}\label{eq:31}
\crazy{a}{pa}^2 (\crazy{a}{qb}  - \crazy{a}{qc})
+ \crazy{a}{pb}^2 (\crazy{a}{qc}  - \crazy{a}{qa})
+ \crazy{a}{pc}^2 (\crazy{a}{qa}  - \crazy{a}{qb})
+ \crazy{a}{pa} \;\crazy{a}{qa} (\crazy{a}{pc} - \crazy{a}{pb})
+ \crazy{a}{pb}\; \crazy{a}{qb} (\crazy{a}{pa} - \crazy{a}{pc})
+ \crazy{a}{pc}\; \crazy{a}{qc} (\crazy{a}{pb} - \crazy{a}{pa})
= 0,
\end{equation}
where $p, q = 1, 2; \;a, b, c = 3, 4, 5; \; p \neq q;
\; a \neq b, c; \; b \neq c.$
\end{lemma}

{\sf Proof.} The proof of (30) and (31) is straightforward and can be obtained by
applying conditions (27) (or (25)) several times. \rule{3mm}{3mm}

\begin{theorem} The components $\crazy{a}{\alpha\beta}$
 of the torsion tensor and their Pfaffian derivatives
$\crazy{a}{\alpha\beta\gamma}$ satisfy the following identities:
\begin{equation}\label{eq:32}
\left\{
\begin{array}{ll}
m_1 = 0, \;\; m_2 = 0, \;\; m_\sigma = 0, \;\;\;\; \sigma = 6, \ldots, n, \\
m_3 =               C  \crazy{a}{34}  + A \crazy{a}{35}, \\
m_4 =  B   \crazy{a}{34}  + A \crazy{a}{45}, \\
m_5 = B \crazy{a}{35}  + C \crazy{a}{45}.
\end{array}
\right.
\end{equation}
where
$$
m_\alpha = (\crazy{a}{24} -  \crazy{a}{15}) \crazy{a}{13\alpha}
+ (\crazy{a}{25} -  \crazy{a}{13}) \crazy{a}{14\alpha}
+ (\crazy{a}{23} -  \crazy{a}{14}) \crazy{a}{15\alpha}
+ (\crazy{a}{15} -  \crazy{a}{14}) \crazy{a}{23\alpha}
+ (\crazy{a}{13} -  \crazy{a}{15}) \crazy{a}{24\alpha}
+ (\crazy{a}{14} -  \crazy{a}{13}) \crazy{a}{25\alpha}
$$
and
$$
A = \left|
\begin{array}{ll}
\crazy{a}{13} & \crazy{a}{14}  \\
\crazy{a}{23} & \crazy{a}{24}
\end{array}
\right|, \;\;
B =
\left|
\begin{array}{ll}
\crazy{a}{14} & \crazy{a}{15}  \\
\crazy{a}{24} & \crazy{a}{25}
\end{array}
\right|, \;\;
C =
\left|
\begin{array}{ll}
\crazy{a}{15} & \crazy{a}{13}  \\
\crazy{a}{25} & \crazy{a}{23}
\end{array}
\right|, \;\; A + B + C = 0.
$$
\end{theorem}

{\sf Proof.} In fact, differentiating (27), we arrive at the
following Pfaffian equation:
\begin{eqnarray}\label{eq:33}
&&(\crazy{a}{24} -  \crazy{a}{25}) \nabla \crazy{a}{13}
+ (\crazy{a}{25} -  \crazy{a}{23}) \nabla \crazy{a}{14}
+ (\crazy{a}{23} -  \crazy{a}{24}) \nabla \crazy{a}{15} \nonumber \\
&+ & (\crazy{a}{15} -  \crazy{a}{14}) \nabla \crazy{a}{23}
+ (\crazy{a}{13} -  \crazy{a}{15}) \nabla \crazy{a}{23}
+ (\crazy{a}{14} -  \crazy{a}{13}) \nabla \crazy{a}{25} = 0,
\end{eqnarray}
where $\nabla \crazy{a}{\alpha\beta} = d\crazy{a}{\alpha\beta}
- \crazy{a}{\alpha\beta} \omega$. Substituting
into (33) the values of
$\nabla \crazy{a}{\alpha\beta} $ taken from equations (4),
 equating to 0 the coefficients in independent
1-forms $\crazy{\omega}{\alpha}$, and simplifying the equations
obtained by means of (30) and (31), we arrive at
conditions (32).  \rule{3mm}{3mm}

\textbf{6.} For the  Goursat web $W (n + 1, n, 1)$ of the second
kind, the second and the third equations of (21) are equivalent
 either to the equation
\begin{equation}\label{eq:34}
  (\crazy{a}{14} - \crazy{a}{13}) \crazy{\omega}{4}
   + (\crazy{a}{15} -  \crazy{a}{13}) \crazy{\omega}{5}= 0
\end{equation}
or to the equation
\begin{equation}\label{eq:35}
  (\crazy{a}{24} - \crazy{a}{23}) \crazy{\omega}{4}
   + (\crazy{a}{25} -  \crazy{a}{23}) \crazy{\omega}{5}= 0.
\end{equation}
By (27),  equations (34) and (35) are equivalent.
Each of these two equations along with equations $\omega_\sigma = 0, \,
\sigma = 6, \ldots , n$,  defines a four-dimensional distribution
$\Delta_4$. The following theorem gives the conditions of
integrability of the distribution $\Delta_4$.

\begin{theorem}
For the  Goursat web $W (n + 1, n, 1)$ of the second
kind, the distribution $\Delta_4$ is integrable if and only if
the following three conditions hold:
\begin{equation}\label{eq:36}
  n_1 = 0, \;\; n_2 = 0, \;\; n_3 = \crazy{a}{13}
  [(\crazy{a}{15} - \crazy{a}{13})  \crazy{a}{34}
  + (\crazy{a}{13} - \crazy{a}{14}) \crazy{a}{35}],
\end{equation}
where
$$
n_h = (\crazy{a}{15} - \crazy{a}{14}) \crazy{a}{13h}
+ (\crazy{a}{13} - \crazy{a}{15}) \crazy{a}{14h}
+  (\crazy{a}{14} - \crazy{a}{13}) \crazy{a}{15h}, \;\; h = 1, 2,
3.
$$
\end{theorem}

{\sf Proof.} Taking the exterior derivative of equation (34),
we obtain the following exterior quadratic equation:
\begin{eqnarray}\label{eq:37}
 && \!\!\!\![\nabla \crazy{a}{14} - \nabla \crazy{a}{13}
 + (\crazy{a}{13} - \crazy{a}{14}) (\crazy{a}{14}\;
\crazy{\omega}{1} + \crazy{a}{24}\; \crazy{\omega}{2}
+ \crazy{a}{34}\;\crazy{\omega}{3})]
 \wedge \crazy{\omega}{4} \nonumber \\
& +& \!\!\!\! [\nabla \crazy{a}{15} - \nabla \crazy{a}{13}
 + (\crazy{a}{13} - \crazy{a}{15}) (\crazy{a}{15} \;
\crazy{\omega}{1} + \crazy{a}{25} \;\crazy{\omega}{2}
+ \crazy{a}{35}\;\crazy{\omega}{3})]
 \wedge \crazy{\omega}{5} = 0.
 \end{eqnarray}
Next, we use (34) to express the form $\crazy{\omega}{4}$ in terms of the form
$\crazy{\omega}{5}$,  substitute its value
into equation (37), and equate
to 0 the coefficients in the independent exterior quadratic products
$\crazy{\omega}{1} \wedge \crazy{\omega}{5}, \,
\crazy{\omega}{2} \wedge \crazy{\omega}{5}$, and
$\crazy{\omega}{3} \wedge \crazy{\omega}{5}$ (there are no other
 exterior quadratic products in the exterior quadratic equation).
 As a result, we obtain conditions (36). \rule{3mm}{3mm}

 \textbf{Remark} Note that the conditions of integrability
 of equation (35) have the form
\begin{equation}\label{eq:38}
  r_1 = 0, \;\; r_2 = 0, \;\; r_3 = \crazy{a}{23}
  [(\crazy{a}{25} - \crazy{a}{23})  \crazy{a}{34}
  + (\crazy{a}{23} - \crazy{a}{24}) \crazy{a}{35}],
\end{equation}
where
$$
r_h = (\crazy{a}{25} - \crazy{a}{24}) \crazy{a}{23h}
+ (\crazy{a}{23} - \crazy{a}{25}) \crazy{a}{24h}
+  (\crazy{a}{24} - \crazy{a}{23}) \crazy{a}{25h}, \;\; h = 1, 2,
3.
$$

However, these conditions are not independent: they follow from
the corresponding conditions (32) and (36). Moreover,
a straightforward calculation shows that for each
$h = 1, 2, 3$, any two of three systems (32), (36), and
(38) imply the third one.

\textbf{7.} The two-dimensional distribution $\Delta_2$ defined
by equations (21) and (22) belongs to the four-dimensional distribution
$\Delta_4$. In general,  $\Delta_2$ is not integrable even if
 $\Delta_4$ is integrable. However, for the  Goursat web $W (n + 1, n, 1)$
 of the second kind,  $\Delta_2$ is  integrable  if
 $\Delta_4$ is integrable.

 \begin{theorem} If a  web $W (n + 1, n, 1)$ is
 the  Goursat web of the second kind, then integrability
 of  the distribution $\Delta_4$ implies  integrability
 of  the distribution $\Delta_2$.
 \end{theorem}

{\sf Proof.} First, we find the conditions of integrability
  of  the distribution $\Delta_2$. Note that
 exterior differentiation of the first equation of (21) and
 equation (22) leads to the identities. Thus, we must take
 the exterior
 derivative of  the second or the third equation of (21).
 However,  we will take the exterior derivative of equation
 (34) which is equivalent to each of them:
\begin{eqnarray}\label{eq:39}
 && \!\!\!\![\nabla \crazy{a}{14} - \nabla \crazy{a}{13}
 + (\crazy{a}{13} - \crazy{a}{14}) (\crazy{a}{14}\;
\crazy{\omega}{1} + \crazy{a}{24}\; \crazy{\omega}{2})]
 \wedge \crazy{\omega}{4} \nonumber \\
 &+& \!\!\!\! [\nabla \crazy{a}{15} - \nabla \crazy{a}{13}
 + (\crazy{a}{13} - \crazy{a}{15}) (\crazy{a}{15} \;
\crazy{\omega}{1} + \crazy{a}{25} \;\crazy{\omega}{2})]
 \wedge \crazy{\omega}{5} = 0.
 \end{eqnarray}
 Using equations (21) and (22), we can express the form
$\crazy{\omega}{2}$ in terms of $\crazy{\omega}{1}$ and
the forms $\crazy{\omega}{3}$ and $\crazy{\omega}{4}$
in terms of the form $\crazy{\omega}{5}$.
As a result, equation (39) will contain only one
independent exterior quadratic product
$\crazy{\omega}{1} \wedge \crazy{\omega}{5}$. Equating
the coefficient in this product to 0, we find that
the condition of integrability is
\begin{equation}\label{eq:40}
  (\crazy{a}{15} - \crazy{a}{14}) (\crazy{a}{131} - \crazy{a}{132})
+ (\crazy{a}{13} - \crazy{a}{15}) (\crazy{a}{141} - \crazy{a}{142})
+ (\crazy{a}{14} - \crazy{a}{13}) (\crazy{a}{151} - \crazy{a}{152})
= 0.
\end{equation}
Now it is easy to see that condition (40) is identically satisfied
if the first two of conditions (36) are satisfied. Thus,
integrability of  $\Delta_4$ implies  integrability
 of  $\Delta_2$. \rule{3mm}{3mm}

\textbf{8.} For the  Goursat web $W (n + 1, n, 1)$ of the second
kind, the first three equations of (23) are equivalent to
 any pair of the following three equations
\begin{equation}\label{eq:41}
  (\crazy{a}{1a} - \crazy{a}{1b}) \crazy{\omega}{1}
   + (\crazy{a}{2a} -  \crazy{a}{2b}) \crazy{\omega}{2}= 0, \;\; a, b
   = 3, 4, 5; \, a \neq b.
\end{equation}
By (27), only two of the three  equations (41) are
independent.
Each  of these equations along with equations $\omega_\sigma = 0, \,
\sigma = 6, \ldots , n$,  defines a four-dimensional distribution
$\Delta_4^\prime$. The following proposition gives the conditions of
integrability of the distribution $\Delta_4^\prime$.

\begin{theorem}
For the  Goursat web $W (n + 1, n, 1)$ of the second
kind, the distribution $\Delta_4^\prime$ is integrable if and only if
the following three conditions hold:
\begin{equation}\label{eq:42}
  s_3 = - C \crazy{a}{34}, \;\; s_4 = B a_{34}, \;\;
  s_5 = C(\crazy{a}{35} - \crazy{a}{45}),
\end{equation}
where
$$
s_h = (\crazy{a}{23} - \crazy{a}{25}) (\crazy{a}{14h}
- \crazy{a}{13h}) + (\crazy{a}{15} -  \crazy{a}{13})
 (\crazy{a}{24h} - \crazy{a}{23h}), \;\; h = 3, 4, 5.
$$
\end{theorem}

{\sf Proof.} We take the exterior derivative of equation (41)
considered for $a = 4, \,b = 3$:
\begin{eqnarray}\label{eq:43}
 && \!\!\!\![\nabla \crazy{a}{14} - \nabla \crazy{a}{13}
 + (\crazy{a}{13} - \crazy{a}{14}) (\crazy{a}{13} \;
\crazy{\omega}{3} + \crazy{a}{14}\; \crazy{\omega}{4}
+ \crazy{a}{15}\; \crazy{\omega}{5})]
 \wedge \crazy{\omega}{1}\nonumber \\
& + & \!\!\!\!  [\nabla \crazy{a}{24} - \nabla \crazy{a}{23}
 + (\crazy{a}{23} - \crazy{a}{24}) (\crazy{a}{23} \;
\crazy{\omega}{3} + \crazy{a}{24}\; \crazy{\omega}{4}
+ \crazy{a}{25}\; \crazy{\omega}{5})]
 \wedge \crazy{\omega}{2} = 0.
\end{eqnarray}
Next, we use (41) to express the form $\crazy{\omega}{2}$ in terms of the form
$\crazy{\omega}{1}$,  substitute its value
into equation (43), and equate
to 0 the coefficients in the independent exterior quadratic products
$\crazy{\omega}{1} \wedge \crazy{\omega}{3}, \,
\crazy{\omega}{1} \wedge \crazy{\omega}{4}$, and
$\crazy{\omega}{1} \wedge \crazy{\omega}{5}$ (there are no other
 exterior quadratic products in the exterior quadratic equation).
 As a result, we obtain conditions (42). \rule{3mm}{3mm}

 \textbf{Remark} Note that the conditions of integrability
 of   equations (41) taken for $a = 5,\, b = 3$
 or $a = 4,\, b = 5$, are obtained from
the corresponding conditions (32) and (42).

\textbf{9}.
The three-dimensional distribution $\Delta_3$ defined
by equations (23) belongs to the four-dimensional distribution
$\Delta_4^\prime$. In general,  $\Delta_3$ is not integrable even if
 $\Delta_4^\prime$ is integrable.
 As opposed to what we had in Theorem 8,
 the following theorem shows that this is also true
  for the  Goursat web $W (n + 1, n, 1)$
 of the second kind.

 \begin{theorem} For  the  Goursat web  $W (n + 1, n, 1)$
 of the second kind, the distribution $\Delta_3$
 is  integrable if and only if the following conditions hold:
\begin{equation}\label{eq:44}
 u_4 = 2 A \crazy{a}{34}, \;\; u_5 = 2 A \crazy{a}{35}, \;\;
 v_4 - u_4 = - A \crazy{a}{34}, \;\; v_5 - u_5 =  A
(\crazy{a}{34} - \crazy{a}{35}),
\end{equation}
where
$$
\begin{array}{ll}
u_k =  (\crazy{a}{23} -  \crazy{a}{24})  \crazy{a}{13k}
- (\crazy{a}{14} -  \crazy{a}{13})  \crazy{a}{23k}, \\
v_k =  (\crazy{a}{23} -  \crazy{a}{24})  \crazy{a}{14k}
- (\crazy{a}{14} -  \crazy{a}{13})  \crazy{a}{24k}, \, k = 4, 5.
\end{array}
$$

 \end{theorem}

{\sf Proof.} In order to find the conditions of integrability
  of  the distribution $\Delta_3$,
 we must take the exterior
 derivatives of only  the first two equations of (23) since
 by (27) the third equation follows from the second one.
 We will take the exterior derivative of the first equation
 of (23) and equation (41) taken for $a = 4,\, b = 3$,  which is equivalent to
 the second equation of (23).

 Exterior differentiation of these two equations gives
the following exterior quadratic equations:
\begin{eqnarray}\label{eq:45}
& & \!\!\!\![\nabla \crazy{a}{13} - \crazy{a}{13} (\crazy{a}{12}\;
\crazy{\omega}{2} + \crazy{a}{13}\; \crazy{\omega}{3} +
\crazy{a}{14}\;
\crazy{\omega}{4} + \crazy{a}{15}\; \crazy{\omega}{5})]
 \wedge \crazy{\omega}{1} \nonumber \\
 &+& \!\!\!\![\nabla \crazy{a}{23} - \crazy{a}{23} (\crazy{a}{21}\;
\crazy{\omega}{1} + \crazy{a}{23}\; \crazy{\omega}{3} +
\crazy{a}{24}\;
\crazy{\omega}{4} + \crazy{a}{25}\; \crazy{\omega}{5})]
 \wedge \crazy{\omega}{2} \nonumber \\
& -& \!\!\!\! (\crazy{a}{13} \;
\crazy{\omega}{1} + \crazy{a}{23} \;\crazy{\omega}{2}
 + \crazy{a}{34} \;\crazy{\omega}{4}  +  + \crazy{a}{35} \;\crazy{\omega}{5})
 \wedge \crazy{\omega}{3} = 0
 \end{eqnarray}
and
\begin{eqnarray}\label{eq:46}
& & \!\!\!\![\nabla \crazy{a}{14} - \nabla \crazy{a}{13}
 + (\crazy{a}{13} - \crazy{a}{14}) (\crazy{a}{12}\;
\crazy{\omega}{2} + \crazy{a}{13}\;
\crazy{\omega}{3} + \crazy{a}{14}\;
\crazy{\omega}{1} + \crazy{a}{24}\; \crazy{\omega}{2})]
 \wedge \crazy{\omega}{1} \nonumber \\
 &+& \!\!\!\! [\nabla \crazy{a}{24} - \nabla \crazy{a}{23}
 + (\crazy{a}{24} - \crazy{a}{23}) (\crazy{a}{21} \;
\crazy{\omega}{1} + \crazy{a}{23} \;\crazy{\omega}{3}
+ \crazy{a}{24} \;\crazy{\omega}{4} + \crazy{a}{25} \;\crazy{\omega}{5})]
 \wedge \crazy{\omega}{2} = 0.
 \end{eqnarray}
 Using equations (23) and (41), we can express the forms
$\crazy{\omega}{2}$ and $\crazy{\omega}{3}$ in terms of
the form $\crazy{\omega}{1}$.
As a result, equations (45) and (46) will contain only two
independent exterior quadratic products
$\crazy{\omega}{1} \wedge \crazy{\omega}{4}$ and
$\crazy{\omega}{1} \wedge \crazy{\omega}{5}$ (there is no
the form $\crazy{\omega}{4} \wedge \crazy{\omega}{5}$ in
(45) and (46)). Equating
the coefficients in these products in equation (45) to 0,
we obtain the first two conditions of (44);
and equating  the coefficients in these products
in equation (46) to 0,
we obtain the last two conditions of (44).
\rule{3mm}{3mm}

It is easy to see that conditions (42) do not imply conditions
(44), i.e., integrability of $\Delta_4^\prime$ does not imply
integrability of $\Delta_3$.

\noindent {\em Author's address}:

\noindent
Vladislav V. Goldberg\\
Department of Mathematics\\
New Jersey Institute of Technology\\
Newark, N.J. 07102, U.S.A.

\vspace*{2mm}

\noindent
 E-mail address: vlgold@m.njit.edu


\begin{thebibliography}{Gou 99}

\bibitem[G 73]{g73} Goldberg, V.~V., {\em $(n+1)$-webs of
                     multidimensional surfaces,}
                     Dokl. Akad. Nauk SSSR \textbf{210} (1973),
                     no. 4, 756--759 (Russian); English translation:
                     Soviet Math. Dokl. \textbf{14} (1973), no. 3,
                     795--799.

\bibitem[G 74]{g74} Goldberg, V.~V., {\em $(n+1)$-webs of
                      multidimensional surfaces,}
                      Bulgar. Akad. Nauk Izv. Mat. Inst. \textbf{15}
                      (1974), 405--424 (Russian).

\bibitem[G 76]{g76} Goldberg, V.~V.,
 {\em Reducible $(n+1)$-webs, group $(n+1)$-webs and
      $(2n+2)$-hedral
      $(n+1)$-webs of multidimensional surfaces,}
      Sibirsk. Mat. Zh. \textbf{17} (1976), no. 1, 44--57
       (Russian); English
  translation: Siberian Math. J. \textbf{17} (1976), no. 1, 34--44.

\bibitem[G 88]{g88} Goldberg, V.~V.,
 {\em Theory of Multicodimensional $(n+1)$-Webs}, Kluwer
      Academic Publishers, Dordrecht-Boston-Tokyo,
      1988, xxii + 466 pp.

\bibitem[Go 99]{go99} Goursat,~E.,
\textit{Sur les \'{e}quations du second ordre \`{a} n
variables, analogues \`{a} l'\'{e}quation de Monge--Amp\`{e}re},
Bull. Soc. Math. France \textbf{27} (1899), 1--34.

\end{thebibliography}
\end{document}